\documentclass[notitlepage]{article}

 \usepackage[matrix,arrow,curve]{xy}
 \usepackage{graphicx}
 \usepackage[english]{babel}
 \usepackage{amsmath}
 \usepackage{amssymb, indentfirst}
 \usepackage{amsthm}
 \usepackage{enumitem}
\begin{document}

\newtheorem{sled}{Corrolary}
\newtheorem{lem}{Lemma}
\newtheorem{zam}{Remark}
\newtheorem{ex}{Example}
\newtheorem{opr}{Definition}
\newtheorem{thm}{Theorem}
\newtheorem{predl}{Proposition}
\setlist[enumerate]{label=\arabic*), font=\normalfont}

\title{TWO-TERM PARTIAL TILTING COMPLEXES OVER BRAUER TREE ALGEBRAS}
\author{Mikhail Antipov and Alexandra Zvonareva}
\date{}
\maketitle

\begin{abstract}

In this paper we describe all indecomposable two-term partial
tilting complexes over a Brauer tree algebra with multiplicity 1
using a criterion for a minimal projective presentation of a module
to be a partial tilting complex. As an application we describe all
two-term tilting complexes over Brauer star algebra and compute
their endomorphism rings.

\end{abstract}

\section{Introduction}

In \cite{RZ} Rouquier and Zimmermann defined a derived Picard group
$\text{TrPic}(A)$ of an algebra $A$, i. e. a group of
autoequivalences of the derived category of $A$, given by
multiplication by a two-sided tilting complex modulo natural
isomorphism. The tensor product of two-sided tilting complexes gives
the multiplication in this group. Despite the fact that for a Brauer
tree algebra with the multiplicity of the exceptional vertex 1
several braid group actions on $\text{TrPic}(A)$ are known
(\cite{RZ}, \cite{IM}), the whole derived Picard group is computed
only for an algebra with two simple modules (\cite{RZ}).

On the other hand Abe and Hoshino showed that over a selfinjective
artin algebra of finite representation type any tilting complex $P$
such that $\text{add}(P) = \text{add}(\nu P)$, where $\nu$ is the
Nakayama functor, can be presented as a product of tilting complexes
of length $\leq 1$ (\cite{AH}). Therefore instead of considering the
derived Picard group we can consider the derived Picard groupoid
corresponding to some class of derived equivalent algebras. The
objects of this groupoid are the algebras from this class and the
morphisms are the derived equivalences given by multiplication by a
two-sided tilting complex modulo natural isomorphism. For example,
one can consider the derived Picard groupoid corresponding to the
class of Brauer tree algebras with fixed number of simple modules
and multiplicity $k$ (the algebras from this class are derived
equivalent and this class is closed under derived equivalence). Then
the result of Abe and Hoshino means that the derived Picard groupoid
corresponding to the class of Brauer tree algebras with fixed number
of simple modules and multiplicity $k$ is generated by one-term and
two-term tilting complexes.

In this paper we give a criterion for a minimal projective
presentation of a module without projective direct summands to be a
partial tilting complex, namely we have the following:

\textbf{Proposition 1}\emph{ Let $A$ be a selfinjective $K$-algebra,
$M$ be a module without projective direct summands and let $T:=  P^0
\overset{f}{\rightarrow} P^1$ be a minimal projective presentation
of module $M.$ Complex $T$ is partial tilting if and only if
$\emph{Hom}_{A}(M,\Omega^2M)=0$ and $\emph{Hom}_{K^b(A)}(T,M)=0.$}

In Proposition 1 module $M$ is considered as a stalk complex
concentrated in degree $0$, complex $T:=  P^0
\overset{f}{\rightarrow} P^1$ is concentrated in degrees $0$ and $1$
accordingly.

Using this proposition we classify all indecomposable two-term
partial tilting complexes over a Brauer tree algebra with
multiplicity 1.

\textbf{Theorem 1}\emph{ Let $A$ be a Brauer tree algebra with
multiplicity 1. A  minimal projective presentation of an
indecomposable non-projective $A$-module $M$ is a partial tilting
complex if and only if  $M$ is not isomorphic to $P/ \emph{soc}(P)$
for any indecomposable projective module $P.$}

Hopefully it will allow us to obtain a full classification of
two-term tilting complexes over Brauer tree algebras. As an
illustration we describe all two-term tilting complexes over Brauer
star algebra and compute their endomorphism rings (for an arbitrary
multiplicity) in sections 5 and 6. Note that the results in sections
5 and 6 partially intersect with \cite{SI1}, \cite{SI2}.

\textbf{Acknowledgement:} We would like to thank Alexander Generalov
for his helpful remarks.

\section{Preliminaries}

Let $K$ be an algebraically closed field, $A$ be a finite
dimensional algebra over $K$. We will denote by $A\text{-}{\rm mod}$
the category of finitely generated left $A$-modules, by $K^b(A)$ --
the bounded homotopy category and by $D^b(A)$ the bounded derived
category of $A\text{-}{\rm mod}.$ The shift functor on the derived
category will be denoted by $[1].$ Let us denote by $A\text{-}{\rm
perf}$ the full subcategory of $D^b(A)$ consisting of perfect
complexes, i.e. of bounded complexes of finitely generated
projective $A$-modules. In the path algebra of a quiver the product
of arrows $\overset{a}{\rightarrow} \overset{b}{\rightarrow}$ will
be denoted by $ab.$ For convenience all algebras are supposed to be
basic.

\begin{opr}
A complex $T \in A\text{-}{\rm perf}$ is called tilting if
\begin{enumerate}
    \item $\emph{Hom}_{D^b(A)}(T,T[i])=0, \mbox{ for } i \neq 0$;
    \item T \mbox{ generates }$A\text{-}{\rm perf}  \mbox{ as a triangulated category.}$
\end{enumerate}

\end{opr}

Tilting complexes were defined by Rickard (\cite{Ri1}) and play an
essential role in the study of the equivalences of derived
categories.

\begin{opr}
A complex $T \in A\text{-}{\rm perf}$ is called partial tilting if
the condition $\emph{1}$ from definition $\emph{1}$ is satisfied.
\end{opr}

\begin{opr}
A tilting complex $T \in A\text{-}{\rm perf}$ is called basic if it
does not contain isomorphic direct summands or equally if
$\emph{End}_{D^b(A)}(T)$ is a basic algebra.

\end{opr}

We will call a (partial) tilting complex a two-term (partial)
tilting complex if it is concentrated in two neighboring degrees.

\begin{opr}
An algebra $A$ is called special biserial (\emph{SB}-algebra), if
$A$ is isomorphic to  $KQ/I$ for some quiver $Q$ and an admissible
ideal of relations $I,$ and the following is satisfied:
\begin{enumerate}
    \item any vertex of $Q$ is the starting point of at most two arrows;
    \item any vertex of $Q$ is the end point of at most two arrows;
    \item if $b$ is an arrow in $Q$ then there is at most one arrow $a$ such that $ab \notin I$;
    \item if $b$ is an arrow in $Q$ then there is at most one arrow $c$ such that $bc \notin I$.
\end{enumerate}

\end{opr}

For an SB-algebra the full classification of indecomposable modules
up to isomorphism is known (\cite{GP}, \cite{WW}).

\begin{opr}

Let $B$ be a symmetric \emph{SB}-algebra over a field $K.$ $A$-cycle
is a maximal ordered set of nonrepeating arrows of $Q$ such that the
product of any two neighboring arrows is not equal to zero.

\end{opr}

Note that the fact that algebra is symmetric means that $A$-cycles
are actually cycles. Also sometimes just a maximal ordered set of
arrows of $Q$ such that the product of any two neighboring arrows is
not equal to zero is called an $A$-cycle (see \cite{AG}). Note also
that in this case $A$-cycles are maximal nonzero paths in $B$.

An important example of an SB-algebra of finite representation type
is a Brauer tree algebra. Also these algebras play an important role
in modular representation theory of finite groups.

\begin{opr}

Let $\Gamma$ be a tree with $n$ edges and an exceptional vertex
which has an assigned multiplicity $k \in \mathbb{N}$. Let us fix a
cyclic ordering of the edges adjacent to each vertex in $\Gamma$ (if
$\Gamma$ is embedded into plane we will assume that the cyclic
ordering is clockwise). In this case $\Gamma$ is called a Brauer
tree of type $(n,k)$.

\end{opr}

For a Brauer tree of type $(n,k)$ one can associate a finite
dimensional algebra $A(n,k)$. Algebra $A(n,k)$ is an algebra with
$n$ simple modules $S_i$ which are in one to one correspondence with
edges $i \in \Gamma$. The two series of composition factors of an
indecomposable projective module $P_i$ (with top $S_i$) are obtained
by going anticlockwise around the $i$-th vertex. We go around the
$i$-th vertex $k$ times if the vertex is exceptional and one time if
it is not. The full description of the Brauer tree algebras in terms
of composition factors is given in \cite{Al}.

Furthermore, Rickard showed that two Brauer tree algebras
corresponding to the trees $\Gamma$ and $\Gamma'$ are derived
equivalent if and only if their types $(n,k)$ and $(n',k')$ coincide
(\cite{Ri2}) and it follows from the results of Gabriel and
Riedtmann that this class is closed under derived equivalence
(\cite{GR}).

\section{Two-term tilting complexes over selfinjective algebras}

Let $A$ be an arbitrary finite dimensional  selfinjective
$K$-algebra.

\begin{lem}

Any two-term complex $T:=  P^0 \overset{f}{\rightarrow} P^1 \in
A\text{-}{\rm perf}$ is isomorphic to a direct sum of the minimal
projective presentation of a module and a stalk complex of
projective module concentrated in degree 0.

\end{lem}

\textbf{Proof} Let us denote by $M$ the cokernel of $f$. The minimal
projective presentation of $M$ is a direct summand of $T$. So $T$ is
a direct sum of the minimal projective presentation of $M$, some
stalk complex $P^0$ concentrated in degree 0, which can be zero and
on which $f$ acts as a zero map, and a complex of the form $P
\overset{\text{id}}{\rightarrow} P,$ which is homotopic to
0.~\hfill\(\Box\)

We will suppose that the minimal projective presentation of a module
is concentrated in degrees 0 and 1 in cohomological notation. For
the sake of simplicity we will consider only minimal projective
presentations of modules without projective summands. Direct
summands corresponding to stalk complexes of projective modules
concentrated in degree 1 will be considered separately in
Proposition 2.

\begin{predl}

Let $A$ be a selfinjective $K$-algebra, $M$ be a module without
projective direct summands and let $T:=  P^0
\overset{f}{\rightarrow} P^1$ be a minimal projective presentation
of module $M.$ Complex $T$ is partial tilting if and only if
$\emph{Hom}_{A}(M,\Omega^2M)=0$ and $\emph{Hom}_{K^b(A)}(T,M)=0.$

\end{predl}

\textbf{Proof} Let $h:P^1 \rightarrow P^0$ be such that $hf=0=fh,$
i.e. $h$ gives a morphism $T \rightarrow T[-1].$

$$ \xymatrix { 0 \ar[r]& \text{Ker}(f) \ar[r]^-{i} & P^0 \ar[r]& P^1 \ar[r]^-{\pi} \ar[ld]& \text{Coker}(f) \ar[r] \ar@{-->}[lld] \ar@{-->}[llld]& 0 \\
0 \ar[r]& \text{Ker}(f) \ar[r] & P^0 \ar[r]& P^1 \ar[r] &
\text{Coker}(f) \ar[r] & 0\\}
 $$

The condition $hf=0$ means that $\text{Im}(f)\subseteq
\text{Ker}(h)$, consequently $h$ goes through $\text{Coker}(f)$,
i.e. there exists $h' \in \text{Hom}_{A}(\text{Coker}(f),P^0)$ such
that $h=h'\pi,$ but $\pi$ is surjective, hence
$\text{Im}(h)=\text{Im}(h').$

The condition $fh=0$ means that $\text{Im}(h')=\text{Im}(h)
\subseteq \text{Ker}(f)$ consequently $h'$ goes through
$\text{Ker}(f),$ i.e. there exists $h''$ such that $h'=ih'',$
$h=ih''\pi.$ Note that since $\pi$ is surjective and $i$ is
injective, $h=0$ if and only if $h''=0.$

Also if there is a nonzero $h'' \in \text{Hom}_A(\text{Coker}(f),
\text{Ker}(f))$ a morphism $h=ih''\pi$ gives a nonzero morphism $T
\rightarrow T[-1].$ So \begin{equation}\label{eqn:1}
\text{Hom}_{D^b(A)}(T,T[-1])=0 \Leftrightarrow
\text{Hom}_{A}(M,\Omega^2M)=0. \tag{$\ast$}\end{equation} Let us now
verify that
\begin{equation}\label{eqn:2}\text{Hom}_{D^b(A)}(T,T[1])=0
\Leftrightarrow \text{Hom}_{K^b(A)}(T,M)=0.
\tag{$\ast\ast$}\end{equation} We have that
$\text{Hom}_{D^b(A)}(T,T[1])=\text{Hom}_{D^b(A)}(T,P_{\bullet})=\text{Hom}_{D^b(A)}(T,M),$
where $P_{\bullet}$ is the projective resolution of $M.$ Since $T$
consists of projective modules,
$\text{Hom}_{D^b(A)}(T,M)=\text{Hom}_{K^b(A)}(T,M).$ \hfill\(\Box\)

\begin{sled}

The projective presentation of a band-module over a symmetric
\emph{SB}-algebra can not be a partial tilting complex.

\end{sled}

\textbf{Proof} In the Auslander-Reiten quiver all band-modules lie
on 1-tubes (\cite{BR}), so $\Omega^2M=M.$ \hfill\(\Box\)

The proof of the next statement is analogous to the proof of
Proposition 1.

\begin{predl} Let $A$ be a selfinjective $K$-algebra, $M$ be a module without
 projective direct summands such that its minimal projective
presentation is a partial tilting complex.

The sum of a stalk complex of projective module $P$ concentrated in
degree 0 and the minimal projective presentation of module $M$ is a
partial tilting complex if and only if
$\emph{Hom}_{A}(M,P)=0=\emph{Hom}_{A}(P,M).$

The sum of a stalk complex of projective module $P$ concentrated in
degree 1 and the minimal projective presentation of module $M$ is a
partial tilting complex if and only if
$\emph{Hom}_{A}(\Omega^2M,P)=0=\emph{Hom}_{A}(P,\Omega^2M).$

\end{predl}

\section{Two-term tilting complexes over Brauer tree algebras with multiplicity 1}

The next remark (\cite{Ha}) plays an important role.

\begin{zam}

Let $A$ be a finite dimensional algebra over a field $K$, let
$\text{\rm proj-}A$ and $\text{\rm inj-}A$ be the categories of
finitely generated projective and injective modules respectively,
$K^b(\text{\rm proj-}A)$, $K^b(\text{\rm inj-}A)$ bounded homotopy
categories, $D$ the duality of the module category with respect to
$K.$ Then the Nakayama functor $\nu$ induces an equivalence of
triangulated categories $K^b(\text{\rm proj-}A) \rightarrow
K^b(\text{\rm inj-}A)$ and there is a natural isomorphism
$D\emph{Hom}(P,-)\rightarrow \emph{Hom}(-,\nu P)$ for $P \in
K^b(\text{\rm proj-}A)$.

\end{zam}

In the case of the symmetric algebra it means that for $T \in
A\text{-}{\rm perf}$ the condition $\text{Hom}_{D^b(A)}(T,T[1])=0$
is satisfied if and only if $\text{Hom}_{D^b(A)}(T,T[-1])=0.$

From now on in this section we will consider only Brauer tree
algebras $A$ corresponding to a Brauer tree $\Gamma$ such that the
multiplicity of the exceptional vertex of $\Gamma$ is 1. Let us fix
an $A$-module $M$ and let us denote by $T:=  P^0
\overset{f}{\rightarrow} P^1$ its minimal projective presentation.

\begin{lem}

Let $M$ be an indecomposable nonprojective $A$-module. The condition
$\emph{Hom}_{A}(P^0,M)=0$ implies $\emph{Hom}_{A}(M,\Omega^2M)=0$
and $\emph{Hom}_{K^b(A)}(T,M)=0.$

\end{lem}

\textbf{Proof} The condition $\text{Hom}_{A}(P^0,M)=0$ obviously
implies $\text{Hom}_{K^b(A)}(T,M)=0.$

Let us show that $\text{Hom}_{A}(P^0,M)=0$ implies
$\text{Hom}_{A}(M,\Omega^2M)=0.$ Since $\text{Hom}_{A}(P^0,M)=0$,
there is no composition factor in $M$ isomorphic to a direct summand
of $\text{top}(P^0)=\text{soc}(P^0).$ The module $\Omega^2M$ is a
submodule of $P^0,$ hence $\text{soc}(\Omega^2M) \subseteq
\text{soc}(P^0)$. For any $h \in \text{Hom}_{A}(M,\Omega^2M)$ we
have that $\text{Im}(h)\cap\text{soc}(\Omega^2M)=0,$ hence $h=0.$
\hfill\(\Box\)

\begin{lem}

Let $M$ be a nonprojective $A$-module such that
$\emph{dim}(\emph{top}(M))=1.$ The minimal projective presentation
of $M$ is a partial tilting complex if and only if $M$ is not
isomorphic to $P/\emph{soc}(P)$ for any indecomposable projective
module $P.$

\end{lem}

\textbf{Proof} The condition $\text{dim}(\text{top}(M))=1$ implies
that $M \simeq P^1/U,$ where $P^1$ is indecomposable.

If $U=\text{soc}(P^1),$ then $P^0 \simeq P^1$ because $A$ is
symmetric. Hance $\Omega^2M$ is a submodule of $P^1,$ hence,
$\text{soc}(\Omega^2M)=\text{soc}(P^1)=\text{top}(P^1)=\text{top}(M),$
which means that $\text{Hom}_{A}(M,\Omega^2M) \neq 0.$ By
\eqref{eqn:1} we get that $\text{Hom}_{D^b(A)}(T,T[-1]) \neq 0.$

Let us assume that $U \neq \text{soc}(P^1).$ We denote by $I$ the
set of indexes corresponding to composition factors of
$\text{top}(U).$ The projective cover of $U$ is isomorphic to
$\bigoplus_{i \in I}Ae_i.$ Since $U \neq \text{soc}(P^1)$, the set
$I$ does not contain the indexes corresponding to $\text{soc}(P^1)$
or to composition factors of $P^1/U$ (over a Brauer tree algebra
with multiplicity 1 all composition factors of an indecomposable
projective module except for the top and the socle are distinct).
Hance $\text{Hom}_{A}(P^0,M)=0.$ By Lemma 2 and Proposition 1 the
minimal projective presentation of $P^1/U$ is a partial tilting
complex. \hfill\(\Box\)

Let us denote by $CF(L)$ the set of the composition factors of
module $L$.

\begin{lem}

For any indecomposable nonprojective $A$-module $M$ such that
$\emph{dim}(\emph{top}(M)) \geq 2$ the condition
$\emph{Hom}_{K^b(A)}(T,M)=0$ is satisfied.

\end{lem}

\textbf{Proof} Note that $\text{dim}(\text{top}(M)) \geq 2$ implies
$CF(\text{top}(P^0))\cap CF(M) \subseteq \text{soc}(M)$. Indeed,
since over a Brauer tree algebra with multiplicity 1 all composition
factors of an indecomposable nonprojective module are distinct,
$CF(\text{top}(P^0))\cap CF(M) \subseteq \text{soc}(M)$.
Consequently, for any morphism $h:P^0\rightarrow M$ the following
holds $\text{Im}(h)\subseteq \text{soc}(M),$ hence $\text{Ker}h
\supseteq \text{rad}(P^0) \supseteq \text{Ker}f,$ hence $h$ goes
through $f$ and $h=0$ in $K^b(A).$\hfill\(\Box\)

Finally we have:

\begin{thm}

A  minimal projective presentation of an indecomposable
non-projective $A$-module $M$ is a partial tilting complex if and
only if  $M$ is not isomorphic to $P/ \emph{soc}(P)$ for any
indecomposable projective module $P.$

\end{thm}

\textbf{Proof} The case $\text{dim}(\text{top}(M)) = 1$ is dealt
with in Lemma 3; in the case $\text{dim}(\text{top}(M)) \geq 2$ the
required result holds because of Lemma 4 Remark 1 and \eqref{eqn:2}.
\hfill\(\Box\)

\section{Two-term tilting complexes over Brauer star algebra}

Let us consider a quiver $Q:$
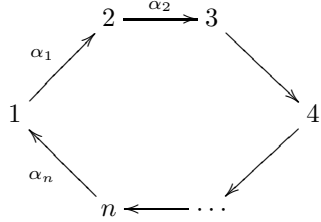
$$
\xymatrix {
& 2 \ar[r]^-{\alpha_2} & 3 \ar[dr] &  \\
    1 \ar[ur]^-{\alpha_1} &  &  & 4 \ar[dl] \\
     & n \ar[ul]^-{\alpha_n} & \cdots \ar[l] & \\
}
$$

The vertices of the quiver are numbered by elements of
$\mathbb{Z}/n\mathbb{Z}.$ Consider the ideal $I$ generated by
relations $$I:=\langle (\alpha_{i} \cdot \alpha_{i+1} \cdot \ldots
\cdot \alpha_{i-1})^k \cdot \alpha_{i}, \mbox{ } i=1,\ldots,n
\rangle.$$ Set $A=kQ/I.$ We denote by $e_i$ the path of length 0
corresponding to the vertex~$i$.

Any indecomposable module over this algebra is uniserial, in
particular any indecomposable module is uniquely determined by the
ordered set of its composition factors. We will denote a module by
the set of the indexes corresponding to its composition factors
ordered from the top to the socle. For example, the simple module
corresponding to the idempotent $e_i$ will be denoted by $(i)$.

In the previous section the description of all two-term partial
tilting complexes in the case $k=1$ was given. Now we will describe
such complexes over a Brauer star algebra for an arbitrary $k.$

\begin{predl}

The minimal projective presentation of an indecomposable $A$-module
is a partial tilting complex if and only if $l(M) < n,$ where $l(M)$
is the length of $M.$

\end{predl}

\textbf{Proof} If $|CF(M)|>n-1$ then both $M$ and $\Omega^2M$
contain all simple modules as composition factors. In particular,
$\text{top}(M)$ is a composition factor of $\Omega^2M$ hence
$\text{Hom}_{A}(M,\Omega^2M) \neq 0.$ If $|CF(M)|<n,$ then in
$\Omega^2M$ there is no composition factor isomorphic to
$\text{top}(M)$ hence $\text{Hom}_{A}(M,\Omega^2M) = 0.$ It is also
clear that $\text{Hom}_{K^b(A)}(T,M)=0,$ since there is no
composition factor isomorphic to $\text{top}(P^0)$ in $M$.
\hfill\(\Box\)

Let us describe all two-term tilting complexes over $A,$
concentrated in degrees 0 and 1. Let there be given two modules
$M=(i, i-1,...,j)$ and $N=(m, m-1,...,l)$ such that the number of
composition factors of $M$ and of $N$ is less then $n.$ Let $T$ be
the minimal projective presentation of $M,$ $T'$ be the minimal
projective presentation of $N.$ Note that $\Omega^2M=(i-1,...,j-1),$
$\Omega^2N=(m-1,...,l-1).$ Let us state when the sum of the minimal
projective presentations of $M$ and $N$ is a partial tilting
complex.

$\text{Hom}_{A}(M,\Omega^2N)=0$ if and only if $i \notin \{m-1,
m-2,...,l-1\}$ or $i \in \{m-1, m-2,...,l-1\},$ but $j \in \{i,
i-1,...,l\}.$

$\text{Hom}_{A}(N,\Omega^2M)=0$ if and only if $m \notin \{i-1,
i-2,...,j-1\}$ or $m \in \{i-1, i-2,...,j-1\},$ but $l \in \{m,
m-1,...,j\}.$

Analysing these conditions we conclude that either the sets $\{i,
i-1,...,j-1\},$ $\{m, m-1,...,l-1\}$ do not intersect or one lies
inside the other.

Now let us figure out when a sum of the minimal projective
presentation of a module $M=(i, i-1,...,j)$ and a stalk complex of a
projective module $P_m=(m, m-1,...,m)$ concentrated in degree 0 is a
partial tilting complex.

$\text{Hom}_{A}(M,P)=0=\text{Hom}_{A}(P,M)$ if and only if $m\notin
\{i, i-1,...,j\}.$

Similarly, a sum of the minimal projective presentation of a module
$M=(i, i-1,...,j)$ and a stalk complex of a projective module
$P_m=(m, m-1,...,m)$ concentrated in degree 1 is a partial tilting
complex if and only if
$\text{Hom}_{A}(\Omega^2M,P)=0=\text{Hom}_{A}(P,\Omega^2M),$ i.e.
$m\notin \{i-1,i-2,...j-1\}.$

Note also that all stalk complexes of projective modules are
concentrated either in degree 0 or in degree 1, since for any two
projective modules $P_m,$ $P_l$ over a Brauer star algebra
$\text{Hom}_{A}(P_m,P_l)\neq 0$.

It is known that in the case of a symmetric algebra of finite
representation type any partial tilting complex with $n$ (where $n$
is the number of isoclasses of simple modules) nonisomorphic direct
summands is tilting (\cite{AH}). Thus to describe all two-term
tilting complexes is the same as to describe all configurations of
$n$ pairwise orthogonal indecomposable complexes, each of which is
either a minimal projective presentation of an indecomposable module
$M$ such that the number of composition factors of $M$ is less then
$n$ or a stalk complex of a projective module concentrated in degree
0 or degree 1, i.e. of $n$ complexes which pairwisely satisfy the
conditions stated before.

We will call an interval a set of vertices of an $n$-gon taken in
order with marked starting point and end point. The covering $S$ of
an $n$-gon by distinguished intervals is the following structure: an
$n$-gon with a partition of its vertices into noncrossing intervals
(we call them outer), each interval can contain from $1$ to $n$
vertices; in each outer interval containing $r$ ($r>1$) vertices
$r-2$ inner intervals are additionally chosen, each of which
contains more that 1 vertex; inner intervals either do not intersect
or lie one inside the other. Also in each outer interval $(i,
i-1,...,j)$ with length greater than 1 we will pick out an interval
of length 1 as follows: it is either a starting point for all outer
intervals or an end point. Note that the covering contains exactly
$n$ intervals. To such a covering $S$ one can assign a two-term
tilting complex $T_{S}$ as follows.

We will consider two cases:

1) To all outer intervals $(i, i-1,...,j) \in S$ of length greater
than 1 an inner interval $(j)$ of length 1 is assigned. Let us
construct a tilting complex as follows: for each interval $(i,
i-1,...,j)$ containing more than 1 vertex take a module $M=(i,
i-1,...,j+1)$, as a direct summand of the tilting complex take the
minimal projective presentation of $M$. For each interval containing
1 vertex take a stalk complex of the projective module corresponding
to this vertex concentrated in degree 0. In this way we get $n$
summands.

2) To all outer intervals $(i, i-1,...,j) \in S$ of length greater
than 1 an inner interval $(i)$ of length 1 is assigned. As before
for each interval $(i, i-1,...,j)$ containing more than 1 vertex
take a module $M=(i, i-1,...,j+1)$, as a direct summand of the
tilting complex take the minimal projective presentation of $M$. For
each interval containing 1 vertex take a stalk complex of the
projective module corresponding to this vertex concentrated in
degree 1. In this way we get $n$ summands.

To the trivial covering, containing only intervals of length $1$,
two tilting complexes: $A$ and $A[-1]$ are assigned.

Based on the previous construction we get the following:

\begin{predl}

Over a Brauer star algebra with $n$ vertices and multiplicity $k$
the set of all basic two-term tilting complexes not isomorphic to
$A$ or $A[-1]$ is in one to one correspondence with the set of all
nontrivial coverings of an $n$-gon by distinguished intervals.

\end{predl}

\section{Endomorphism rings}

Let us construct the endomorphism ring of a two-term tilting complex
over a Brauer star algebra with $n$ vertices and multiplicity $k$,
i.e. the endomorphism rings of a tilting complex corresponding to
the covering $S$ of an $n$-gon. It is well known that it is
isomorphic to a Brauer tree algebra corresponding to some Brauer
tree $\Gamma$ with multiplicity $k.$ For this purpose we first
compute the Cartan matrix of the algebra $\text{End}_{K^b(A)}(T_S)$.
It will tell us which edges of $\Gamma$ are incident to one vertex.
After that we will only have to establish the cyclic ordering of the
edges incident to each vertex of $\Gamma$. It is easy to compute the
Cartan matrix of $\text{End}_{K^b(A)}(T_S)$ using the well known
formulae by Happel \cite{Ha2}: let $Q=(Q^r)_{r \in \mathbb{Z}},
R=(R^s)_{s \in \mathbb{Z}} \in A\text{-}{\rm perf}$, then
$$\sum_i (-1)^i {\rm dim}_K {\rm Hom}_{K^b(A)}(Q,R[i])=\sum_{r,s} (-1)^{r-s}{\rm dim}_K {\rm Hom}_{A}(Q^r,R^s).$$
Note that if ${\rm Hom}_{K^b(A)}(Q,R[i])=0, i \neq 0$ (for example,
in the case when $Q$ and $R$ are summands of a tilting complex) then
the left hand side of the formulae becomes ${\rm dim}_K {\rm
Hom}_{K^b(A)}(Q,R).$

As before we will consider two cases:

1) To all outer intervals $(i, i-1,...,j) \in S$ of length greater
than 1 an inner interval $(j)$ of length 1 is assigned, i.e. all
stalk complexes of projective modules which are direct summands of
$T_S$ are concentrated in degree 0. Let $(i, i-1,...,j)$, $(t,
t-1,...,l) \in S$ be two arbitrary intervals of the covering $S$ of
length greater then 1. And let $(m), (r) \in S$ be intervals of
length 1. It is easy to see that
$$\text{dim}(\text{Hom}_{K^b(A)}(P_j\rightarrow P_i,P_l\rightarrow P_t))=$$ $$=\left\{%
\begin{array}{ll}
    0, & \hbox{if } \{i, i-1,...,j\}\cap\{t, t-1,...,l\}=\varnothing;\\
    0, & \hbox{if } \{i, i-1,...,j\} \subset \{t, t-1,...,l\}, i\neq t, j\neq l; \\
    1, & \hbox{if } \{i, i-1,...,j\} \subset \{t, t-1,...,l\}, i = t, j\neq l; \\
    1, & \hbox{if } \{i, i-1,...,j\} \subset \{t, t-1,...,l\}, i\neq t, j= l; \\
    2, & \hbox{if } \{i, i-1,...,j\} = \{t, t-1,...,l\}. \\
\end{array}%
\right.    $$

$$\text{dim}(\text{Hom}_{K^b(A)}(P_m,P_j\rightarrow P_i))=\left\{%
\begin{array}{ll}
    1, & \hbox{if } m=j; \\
    0, & \hbox{if } m \neq j. \\
\end{array}%
\right. $$
$$\text{dim}(\text{Hom}_{K^b(A)}(P_j\rightarrow P_i,P_m))=\left\{%
\begin{array}{ll}
    1, & \hbox{if } m=j; \\
    0, & \hbox{if } m \neq j. \\
\end{array}%
\right.    $$

$$\text{dim}(\text{Hom}_{K^b(A)}(P_m,P_r))=\left\{%
\begin{array}{ll}
    k, \mbox{ for } m\neq r; \\
    k+1, \mbox{ otherwise}. \\
\end{array}%
\right.$$

These data give us the partition of the vertices of
$\text{End}_{K^b(A)}(T_S)$ into $A$-cycles or equally which edges of
the Brauer tree of algebra $\text{End}_{K^b(A)}(T_S)$ are incident
to the same vertex (we will identify the edges of the Brauer tree of
$\text{End}_{K^b(A)}(T_S)$ and the indecomposable summand of $T_S$
corresponding to them). Now we have to find out the cyclic ordering
of the edges incident to one vertex and which vertex is exceptional.
Note that if we arrange the vertices of the $A$-cycle of length $r$
in such a manner that successive composition of $kr$ morphisms (in
the case of the exceptional vertex) or of $r$ morphisms (in the case
of a nonexceptional vertex) between them is not homotopic to zero
then this arrangement will give us the desired cyclic order.

In the case when all stalk complexes of projective modules are
concentrated in degree 0 in the algebra $\text{End}_{K^b(A)}(T_S)$,
the following types of $A$-cycles can occur: a) the $A$-cycle of
projective modules; b) an $A$-cycle containing an indecomposable
stalk complex of a projective module $P$ concentrated in degree 0
and two-term complexes having $P$ as a 0-component; c) an $A$-cycle
containing two-term complexes with the same 0-components; d) an
$A$-cycle containing two-term complexes with the same components in
degree 1.

For convenience let us use the following notation: a homomorphism
$P_l\rightarrow P_m$ induced by multiplication on the right by
$\alpha_{l}\alpha_{l+1}...\alpha_{m-1}$ will be denoted by
$\alpha_{l,m-1}.$

a) Let $(m_1), (m_2),..., (m_r) \in S$ where the set
$\{m_1,m_2,...,m_r\}$ is ordered according to the cyclic ordering of
the edges in the Brauer star, $r$ is maximal. It is clear that the
following diagram of chain maps holds:
$$\xymatrix {
... \ar[r] & 0 \ar[r] \ar[d]&P_{m_1} \ar[r] \ar[d]^{\alpha_{m_1,m_2-1}}& 0 \ar[d] \ar[r]& ...\\
... \ar[r] & 0 \ar[r] \ar[d]&P_{m_2} \ar[r] \ar[d]& 0 \ar[d] \ar[r]& ...\\
... \ar[r] &... \ar[r] \ar[d]& ... \ar[d] \ar[r]& ... \ar[r] \ar[d]& ...\\
... \ar[r] & 0 \ar[r] \ar[d]&P_{m_r} \ar[r] \ar[d]^{\alpha_{m_r,m_1-1}}& 0 \ar[d] \ar[r]& ...\\
... \ar[r] & 0 \ar[r]&P_{m_1} \ar[r] & 0 \ar[r]& ... \\
}$$ The successive composition of any $kr$ morphisms is not
homotopic to 0. So the edges of $\text{End}_{K^b(A)}(T_S)$
corresponding to stalk complexes of projective modules have a common
vertex and the cyclic ordering in $\text{End}_{K^b(A)}(T_S)$ is
induced by the cyclic ordering in the Brauer star. The vertex of
$\text{End}_{K^b(A)}(T_S)$ corresponding to this cycle is
exceptional.

b) Let $(m_1, m_1 - 1,...,j ), (m_2, m_2 - 1,...,j),..., (m_r, m_r -
1,...,j), (j) \in S,$ where the set $\{j,m_1,m_2,...,m_r\}$ is
ordered according to the cyclic ordering of the edges in the Brauer
star, $r$ is maximal. Let us consider the following diagram of chain
maps:
$$\xymatrix {
... \ar[r] & 0 \ar[r] \ar[d]&P_{j} \ar[r] \ar[d]^{(\alpha_{j,j-1})^k}& 0 \ar[d] \ar[r] & 0 \ar[d] \ar[r]& ...\\
... \ar[r] & 0 \ar[r] \ar[d]&P_{j} \ar[r] \ar[d]^{1} & P_{m_1} \ar[d]^{\alpha_{m_1,m_2-1}} \ar[r] & 0 \ar[d] \ar[r]& ...\\
... \ar[r] & 0 \ar[r] \ar[d]&P_{j} \ar[r] \ar[d] & P_{m_2} \ar[d] \ar[r] & 0 \ar[d] \ar[r]& ...\\
... \ar[r] &... \ar[r] \ar[d]& ... \ar[d] \ar[r]& ... \ar[r] \ar[d] & ... \ar[d] \ar[r] & ...\\
... \ar[r] & 0 \ar[r] \ar[d]&P_{j} \ar[r] \ar[d]^{1}& P_{m_r} \ar[d] \ar[r]&0 \ar[d] \ar[r]& ...\\
... \ar[r] & 0 \ar[r]&P_{j} \ar[r] & 0 \ar[r]& 0  \ar[r]& ... \\
}$$ The successive composition of any $r+1$ morphisms is not
homotopic to 0. That means that the edges of
$\text{End}_{K^b(A)}(T_S)$ corresponding to this $A$-cycle are
ordered in the following way: $\{P_j, P_j\rightarrow
P_{m_1},P_j\rightarrow P_{m_2},...,P_j\rightarrow P_{m_r}\}.$

c) Similarly, if $(m_1, m_1 - 1,...,j ), (m_2, m_2 - 1,...,j),...,
(m_r, m_r - 1,...,j) \in S$ is the set of intervals corresponding to
some $A$-cycle in $\text{End}_{K^b(A)}(T_S),$ where the set
$\{m_1,m_2,...,m_r\}$ is ordered according to the cyclic ordering of
the edges in the Brauer star, $r$ is maximal, then the edges of
$\text{End}_{K^b(A)}(T_S)$ corresponding to this $A$-cycle are
ordered in the following way: $\{ P_j\rightarrow
P_{m_1},P_j\rightarrow P_{m_2},...,P_j\rightarrow P_{m_r}\}.$

d)Let us now consider an $A$-cycle containing the summand with the
same components in degree 1. Let $(j, j-1,...,m_1), (j,
j-1,...,m_2),..., (j,j-1,...,m_r) \in S,$ where the set
$\{m_1,m_2,...,m_r\}$ is ordered according to the cyclic ordering of
the edges in the Brauer star, $r$ is maximal. Then the following
diagram of chain maps holds:
$$\xymatrix {
... \ar[r] & 0 \ar[r] \ar[d]&P_{m_1} \ar[r] \ar[d]^{\alpha_{m_1,m_2-1}}& P_j \ar[d]^{1} \ar[r] & 0 \ar[d] \ar[r]& ...\\
... \ar[r] & 0 \ar[r] \ar[d]&P_{m_2} \ar[r] \ar[d] & P_{j} \ar[d]^{1} \ar[r] & 0 \ar[d] \ar[r]& ...\\
... \ar[r] &... \ar[r] \ar[d]& ... \ar[d] \ar[r]& ... \ar[r] \ar[d] & ... \ar[d] \ar[r] & ...\\
... \ar[r] & 0 \ar[r] \ar[d]&P_{m_{r-1}} \ar[r] \ar[d]^{\alpha_{m_{r-1},m_r-1}}& P_{j} \ar[d]^{1} \ar[r]&0 \ar[d] \ar[r]& ...\\
... \ar[r] & 0 \ar[r] \ar[d]&P_{m_r} \ar[r] \ar[d]^{0} & P_j \ar[r] \ar[d]^{(\alpha_{j,j-1})^k}& 0  \ar[r] \ar[d]& ... \\
... \ar[r] & 0 \ar[r] &P_{m_1} \ar[r] & P_j \ar[r] & 0 \ar[r]& ...\\
}$$ The successive composition of any $r$ morphisms is not homotopic
to 0. This means that the edges of $\text{End}_{K^b(A)}(T_S)$
corresponding to this $A$-cycle are ordered in the following way:
$\{P_{m_1}\rightarrow P_{j},P_{m_2}\rightarrow
P_{j},...,P_{m_r}\rightarrow P_{j}\}.$

This completes the first case. Since for each of the 4 cases of
$A$-cycles we have described the cyclic ordering of vertices, which
is naturally induced by the cyclic ordering of the vertices in the
Brauer star algebra.

2) Let us consider the second case. To all outer intervals $(i,
i-1,...,j) \in S$ of length greater than 1 an inner interval $(i)$
of length 1 is assigned, i.e. all stalk complexes of projective
modules which are direct summand of $T_S$ are concentrated in degree
1. Let $(i, i-1,...,j)$, $(t, t-1,...,l) \in S$ be two arbitrary
intervals of length greater than 1. And let $(m), (r) \in S$ be
intervals of length 1. It is easy to see that
$$\text{dim}(\text{Hom}_{K^b(A)}(P_j\rightarrow P_i,P_l\rightarrow P_t))=$$ $$=\left\{%
\begin{array}{ll}
    0, & \hbox{if } \{i, i-1,...,j\}\cap\{t, t-1,...,l\}=\varnothing;\\
    0, & \hbox{if } \{i, i-1,...,j\} \subset \{t, t-1,...,l\}, i\neq t, j\neq l; \\
    1, & \hbox{if } \{i, i-1,...,j\} \subset \{t, t-1,...,l\}, i = t, j\neq l; \\
    1, & \hbox{if } \{i, i-1,...,j\} \subset \{t, t-1,...,l\}, i\neq t, j= l; \\
    2, & \hbox{if } \{i, i-1,...,j\} = \{t, t-1,...,l\}. \\
\end{array}%
\right.    $$

$$\text{dim}(\text{Hom}_{K^b(A)}(P_m,P_j\rightarrow P_i))=\left\{%
\begin{array}{ll}
    1, & \hbox{if } m=i; \\
    0, & \hbox{if } m \neq i. \\
\end{array}%
\right.    $$
$$\text{dim}(\text{Hom}_{K^b(A)}(P_j\rightarrow P_i,P_m))=\left\{%
\begin{array}{ll}
    1, & \hbox{if } m=i; \\
    0, & \hbox{if } m \neq i. \\
\end{array}%
\right.    $$

$$\text{dim}(\text{Hom}_{K^b(A)}(P_m,P_r))=\left\{%
\begin{array}{ll}
    k, \mbox{ for } m\neq r; \\
    k+1, \mbox{ otherwise}. \\
\end{array}%
\right.$$

As in the previous case, the exceptional vertex corresponds to the
cycle of stalk complexes of projective modules (this time they are
concentrated in degree 1). All $A$-cycles can be divided into 4
types. For 3 of them (namely a, c, d) we already know the cyclic
ordering. The remaining case is:

e) Let $(j, j-1,...,m_1), (j, j-1,...,m_2),..., (j,j-1,...,m_r), (j)
\in S,$ where the set $\{j,m_1,m_2,...,m_r\}$ is ordered according
to the cyclic ordering of the edges in the Brauer star, $r$ is
maximal. Let us consider the following diagram of chain maps:

$$\xymatrix {
... \ar[r] & 0 \ar[r] \ar[d]&P_{m_1} \ar[r] \ar[d]^{\alpha_{m_1,m_2-1}}& P_j \ar[d]^{1} \ar[r] & 0 \ar[d] \ar[r]& ...\\
... \ar[r] & 0 \ar[r] \ar[d]&P_{m_2} \ar[r] \ar[d] & P_{j} \ar[d]^{1} \ar[r] & 0 \ar[d] \ar[r]& ...\\
... \ar[r] &... \ar[r] \ar[d]& ... \ar[d] \ar[r]& ... \ar[r] \ar[d] & ... \ar[d] \ar[r] & ...\\
... \ar[r] & 0 \ar[r] \ar[d]&P_{m_{r-1}} \ar[r] \ar[d]^{\alpha_{m_{r-1},m_r-1}}& P_{j} \ar[d]^{1} \ar[r]&0 \ar[d] \ar[r]& ...\\
... \ar[r] & 0 \ar[r] \ar[d]&P_{m_r} \ar[r] \ar[d] & P_j \ar[r] \ar[d]^{(\alpha_{j,j-1})^k}& 0  \ar[r] \ar[d]& ... \\
... \ar[r] & 0 \ar[r] \ar[d]&0 \ar[r] \ar[d] & P_j \ar[r] \ar[d]^{1}& 0  \ar[r] \ar[d]& ... \\
... \ar[r] & 0 \ar[r] &P_{m_1} \ar[r] & P_j \ar[r] & 0 \ar[r]& ...\\
}$$ The successive composition of any $r+1$ morphisms is not
homotopic to 0. This means that the edges of
$\text{End}_{K^b(A)}(T_S)$ corresponding to this $A$-cycle are
ordered in the following way: $\{P_{m_1}\rightarrow
P_{j},P_{m_2}\rightarrow P_{j},...,P_{m_r}\rightarrow P_{j}, P_j\}.$
\hfill\(\Box\)

The following is clear from the description of endomorphism rings of
two-term tilting complexes.

\begin{zam}

A two-term tilting complex $T_S$ over a Brauer star algebra with $n$
edges and multiplicity 1, which is not isomorphic to $A$ or $A[-1]$,
gives a derived autoequivalence if and only if the covering $S$ of
the $n$-gon has the following form:
$$(j,j-1,...,j+1), (j,j-1,...,j+2),...,(j,j-1),(j) j=1,...,n$$ or $$(j-1,j-2,...,j), (j-2,j-3,...,j),...,(j+1,j),(j), j=1,...,n.$$

The subgroup of the derived Picard group generated by these
autoeqivalences was studied in \emph{\cite{IM}}.

In the case $k\neq 1$ a two-term tilting complex $T_S$ gives a
derived autoequivalence if and only if the covering $S$ is trivial.

\end{zam}

Let us consider an example of a two-term tilting complex and compute
its endomorphism ring.

\begin{ex}

Let $k=1,n=4.$ And let $S=(1,2,3,4), (2,3,4), (2,3), (1).$
\end{ex}

Then $T_S$ consists of the following direct summand: $P_4\rightarrow
P_1, P_4\rightarrow P_2, P_3\rightarrow P_2$ and $P_1,$ concentrated
in degree 1. Let us denote the vertices of
$\text{End}_{K^b(A)}(T_S)$ as follows: $a$ is a vertex corresponding
to $P_4\rightarrow P_1,$ $b$ to $P_4\rightarrow P_2,$ $c$ to
$P_3\rightarrow P_2,$ $d$ to $P_1.$ Then the quiver of
$\text{End}_{K^b(A)}(T_S)$ is of the following form:

$$d \leftrightarrows a \leftrightarrows b\leftrightarrows c,$$
and the Brauer graph is a string: $\bullet \hrulefill \bullet
\hrulefill \bullet \hrulefill \bullet \hrulefill \bullet.$

\begin{predl}

For any algebra $B$ corresponding to the Brauer tree $\Gamma$ with
$n$ edges and multiplicity $k$ there is a two-term tilting complex
$T_S$ over $A$ such that $B\simeq \emph{End}_{K^b(A)}(T_S).$

\end{predl}

\textbf{Proof} Let us assume that the root of $\Gamma$ is chosen in
the exceptional vertex, and that $\Gamma$ is embedded in the plane
in such a manner that all nonroot vertices are situated on the plane
lower than the root according to their level (the further from the
root, the lower, all vertices of the same level lie on a horizontal
line). The edges around vertices are ordered clockwise.

Let us number the edges of the tree $\Gamma$ as follows: put 1 on
the right-hand edge incident to the root, on the next edge incident
to the root according to the order put $1+ k_1+1,$ where $k_1$ is
the number of successors of the nonroot end of the edge with label
1. Let the $(i-1)$-st edge incident to the root be labelled with $m$
and let the nonroot vertex incident to the edge with label $m$ have
$k_{m}$ successors, then put on the $i$-th edge incident to the root
label $m+ k_{m} + 1.$

Further on let us put the labels as follows: consider a vertex of an
odd level (a vertex which can be connected to the root by a path of
odd length), let the edge connecting it to the vertex of a higher
level be labelled with $j.$ Put $j+1+ k_1$ on the right-hand edge
incident to this vertex, where $k_{1}$ is the number of successors
of the other end of this edge. Put $j+1 + k_{1} + k_2 +1$ on the
next edge incident to this vertex, where $k_{2}$ is the number of
successors of the other end of this edge. Further on let us put the
labels by induction: let the $(i-1)$-st edge incident to the fixed
vertex be labelled with $m,$ and let the lower end of the next edge
have $k_{m}$ successors, put $m+ k_{m}+ 1$ on the $i$-th edge
incident to this vertex.

Consider a vertex of an even level, let the edge connecting it to
the vertex of a higher level be labelled with $t$ and let the edge
connecting the other end of the edge labelled with $t$ and the
vertex of a higher level be labelled with $j.$ Put $j+1$ on the
right-hand edge incident to this vertex. Put $j+1 + k_{j+1}+1$ on
the next edge incident to this vertex, where $k_{j+1}$ is the number
of successors of the other end of the edge labelled with $j+1$. Let
the $(i-1)$-st edge incident to the fixed vertex be labelled with
$m,$ and let the lower end of the $(i-1)$-st edge incident to the
fixed vertex have $k_{m}$ successors, put $m+ k_{m}+ 1$ on the
$i$-th edge incident to this vertex.

Let us construct a tilting complex over algebra $A$ using a labelled
tree $\Gamma$. Assume that the root of $\Gamma$ has $l$ children and
there are labels $\{n_1, n_2,...,n_l\}$ on the edges incident to the
root. Take stalk complexes of projective modules
$P_{n_1},P_{n_2}...,P_{n_l}$ concentrated in degree 0 as summands of
the tilting complex. Let us consider a vertex of an odd level.
Assume that the edge connecting it to a vertex of a higher level is
labelled by $j,$ the other edges incident to this vertex have labels
$j_1, j_2,...j_h,$ where $h$ is the number of children of this
vertex. In the tilting complex the following direct summands will
correspond to these edges: $P_j\rightarrow P_{j_1}, P_j\rightarrow
P_{j_2},...,P_j\rightarrow P_{j_h}.$

Let us consider a vertex of an even level. Assume that the edge
connecting it to a vertex of a higher level is labelled by $g,$ the
other edges incident to this vertex have labels $g_1, g_2,...,g_d,$
where $d$ is the number of children of this vertex. In the tilting
complex the following direct summands will correspond to these
edges: $P_{g_1}\rightarrow P_g, P_{g_2}\rightarrow
P_g,...,P_{g_d}\rightarrow P_g.$ It is clear that we have the
desired number of summands. Because of the construction this complex
is tilting and the Brauer tree corresponding to its endomorphism
ring is $\Gamma.$

Similarly, we could construct a tilting complex with all the stalk
complexes of projective modules concentrated in degree 1.
\hfill\(\Box\)

\end{document}